\documentclass[12pt]{article}

\usepackage[utf8]{inputenc}

\usepackage{amsfonts}
\usepackage{amssymb}
\usepackage{amsmath}
\usepackage{amsthm}

\def \le {\leqslant}
\def \ge {\geqslant}

\begin{document}
\begin{Large}
\centerline{{\bf Diophantine exponents for systems of linear forms}}
\centerline{{\bf in two variables}}
\vskip+1cm
 \centerline{ by {\bf Nikolay G. Moshchevitin}\footnote{ Research is supported by RFBR grant No.12-01-00681-a and by the grant of Russian Government, project 11.
G34.31.0053.
  }  
  }
\end{Large}

 \vskip+2cm
 \centerline{\bf Abstract}

\begin{small}
 We improve  on Jarn\'{\i}k's  
inequality between
uniform  Diophantine exponent $\alpha $ and
ordinary Diophantine exponent $\beta$
for a system of $ n\ge 2$ real linear forms in two integer variables.
Jarn\'{\i}k (1949, 1954) proved that $\beta \ge \alpha (\alpha -1)$.
In the present paper we give a better bound in the case $\alpha >1$.
We prove that
$$\beta \ge 
\begin{cases}
\frac{1}{2}\left(\alpha^2-\alpha+1+\sqrt{
(\alpha^2-\alpha+1)^2 +4\alpha^2(\alpha-1)}\right)\,\,\,\,\text{if}\,\,\,\,
1\le \alpha \le 2\cr
\frac{1}{2}\left(
\alpha^2-1+\sqrt{(\alpha^2-1)^2+4\alpha (\alpha-1)}\right)
\,\,\,\,\text{if}\,\,\,\,
\alpha \ge 2
\end{cases}.
$$

\end{small}
\vskip+0.5cm


 \vskip+2cm

{ \bf 1. Jarn\'{\i}k's  theorem.}\,\,\,
In this paper
$$\Theta
=\left(
\begin{array}{ccc}
\theta_1^1&\cdots&\theta_1^m\cr
\cdots &\cdots &\cdots \cr
\theta_n^1&\cdots&\theta_n^m
\end{array}
\right)
$$
  stands for a $m\times n$ real martix and  ${\bf x}=(x_1,...,x_m) \in \mathbb{Z}^m$ is an integer vector.
	In the sequel
$|{\bf x}|^{\rm sup}$ means the sup-norm of a vector ${\bf x}$:
$$
|{\bf x}|^{\rm sup}
=\max_{1\le i \le m}|x_i|.
$$
Consider the function
$$
\psi^{\rm sup}_\Theta (t)= \min_{{\bf x} \in \mathbb{Z}^m:\,0<|{\bf x}|^{\rm sup}\le t} \,\,\,
\max_{1\le j\le n}
||\theta_1^{j} x_1+...+\theta_m^jx_j||.
$$
In this paper we suppose that  for every $t\ge 1$ one has
\begin{equation}\label{poss}
\psi^{\rm sup}_\Theta (t)
>0,\,\,\,\,
\forall\, t\ge 1.
\end{equation}
This is a natural condition on the matrix $\Theta$.

The
{\it uniform} Diophantine exponent $\alpha (\Theta)$ is defined as
follows:
\begin{equation}\label{aa}
\alpha(\Theta) = \sup
\{\gamma >0:\,\,\,
\limsup_{t\to +\infty} t^\gamma \psi^{\rm sup}_\Theta (t) <+\infty \},
\end{equation}
From the Minkowski convex body theorem it follows that
\begin{equation}\label{miiha}
\alpha (\Theta)\ge \frac{m}{n}.
\end{equation}
In addition it is a well known fact that in the case 
$ m=1$ one has
$$
\alpha (\Theta) \le 1.
$$

 The
{\it ordinary} Diophantine exponent $\beta (\Theta)$ is defined as
follows:
\begin{equation}\label{bb}
\beta (\Theta) =\sup
\{\gamma >0:\,\,\, \liminf_{t\to +\infty} t^\gamma \psi^{\rm sup}_\Theta (t)
<+\infty\} .
\end{equation}
Obviously
\begin{equation}\label{ab}
\beta (\Theta) \ge \alpha (\Theta).
\end{equation}
This inequality may be considered as a lower bound for $\beta (\Theta)$
in terms of $\alpha(\Theta)$.
V.  Jarn\'{\i}k improved on the trivial bound (\ref{ab}) in several
papers. Probably his first paper dealing with this topic is the paper
\cite{Sze} published in
'Acta Scientarium Mathematicum Szeged' in 1949.
Here we formulate a general result by  Jarn\'{\i}k from \cite{Cze}.

{\bf Theorem A.} (V. Jarn\'{\i}k \cite{Cze})
 \,\,{\it
Suppose that $\Theta$ satisfies (\ref{poss}).
Then

{\rm (i)}  if $m=1$ and $\Theta$ consists of at least two numbers
$\theta_i^1,\theta_k^1$ linearly independent over $\mathbb{Z}$ together with 1,
then
\begin{equation}\label{sim}
\beta (\Theta) \ge \alpha (\Theta)\cdot \frac{\alpha(\Theta)}{1-\alpha(\Theta)}
;
\end{equation}

{\rm (ii)}  if $m=2$ then
\begin{equation}\label{dva}
\beta (\Theta) \ge \alpha (\Theta)\cdot (\alpha (\Theta)-1);
\end{equation}

{\rm(iii)}  in the case $m\ge 3, n\ge 1$ under the additional condition
$\alpha (\Theta) \ge (5m^2)^{m-1}$ one has

\begin{equation}\label{arb}
\beta (\Theta) \ge \alpha (\Theta)\cdot (\alpha(\Theta)^{1/(m-1)}-3).
\end{equation}

}

In the cases $ m=1,n=2$ and $m=2,n=2$ the inequalities of Theorem A 
are the best possible.
In \cite{Lo} M. Laurent proved a general result (so-called 'four exponents theorem')
from which he deduced the following theorem as a corollary.

{\bf Theorem B.} (M. Laurent \cite{Lo})\,\,
{\it

{\rm (i)} 
Suppose that
$\beta \ge \alpha \cdot \frac{\alpha}{1-\alpha},  \frac{1}{2}\le\alpha\le 1$. Then there exists
$\Theta = \binom{\theta_1^1}{\theta_2^1}$ such that
the numbers $1, \theta_1^1, \theta_2^1$ are linearly independent over
$\mathbb{Z}$ and
$\alpha (\Theta) = \alpha, \beta (\Theta) = \beta$.

{\rm (ii)} 
Suppose that
$\beta \ge \alpha \cdot ({\alpha} -{1}),  \alpha\ge 2$. Then there exists
$\Theta = ({\theta_1^1},{\theta_2^1})$ such that
the numbers 
$1, \theta_1^1, \theta_2^1$ are linearly independent over
$\mathbb{Z}$ and
$\alpha (\Theta) = \alpha, \beta (\Theta) = \beta
$.

}

In a recent paper \cite{SS} W.M. Schmidt and L. Summerer
developed a new powerful method of analysis of the successive minima of one-parameter
families of lattices.
This method
enabled them to improve the inequalities (\ref{sim},\ref{arb}) of Theorem A in the  cases $m=1$ amd $n=1$. As a corollary they obtained the following result.

{\bf Theorem C.} (W.M. Schmidt, L. Summerer \cite{SS})\,\,
{\it

{\rm (i)}
Suppose that $m=1, n \ge 2$ and the matrix $\Theta$ consists of numbers
$\theta_1^1,...,\theta_n^1$ linearly independent over $\mathbb{Z}$ together with 1.
Then
 \begin{equation}\label{o1}
\beta (\Theta)  \ge \alpha(\Theta ) \cdot
\frac{\alpha (\Theta)+ n-2 }{(n-1)(1-\alpha (\Theta))}
.
\end{equation}

{\rm (ii)}
Suppose that $n=1, m \ge 2$ and the matrix $\Theta$ consists of numbers
$\theta_1^1,...,\theta_1^m$ linearly independent over $\mathbb{Z}$ together with 1.
Then
\begin{equation}\label{o2}
\beta 
(\Theta)
\ge \alpha (\Theta)\cdot \frac{(m-1)(\alpha(\Theta)-1)}{1+(m-2)\alpha(\Theta)}.
\end{equation}
}

The proof of the main result from \cite{SS}
relies
on K. Mahler's theory of preudocompaund bodies \cite{M}
and deals with 
difficult analysis of special piecewise linear functions.
An alternative easy geometric proof was given by O. German and N. Moshchevitin in \cite{GM}. The inequalities  
(\ref{o1},\ref{o2})
follow from the main result of
\cite{SS} and  transference inequalities by Y. Bugeaud and M. Laurent \cite{BL}.
Here we should note that the method developed by W.M. Schmidt and
L. Summerer in \cite{SS} cannot be  directly applied to the case $m>1,n>1$,
by some geometric reasons. 

One can easily see that in the cases $n=2$ and $m=2$ inequalities (\ref{o1})
and (\ref{o2}) turn into (\ref{sim}) and (\ref{dva}) respectively.

In the case $ m=1, n= 3$  the best known inequality is due to N. Moshchevitin.

{\bf Theorem D.} (N. Moshchevitin \cite{mCze})\,\,{\it
Suppose that $m=1,n=3$
and the collection $   \theta_1^1, \theta_2^1 , \theta_3^1
$ consists of numbers which, together with 1, are linearly
independent over $\mathbb{Z}$. Then
\begin{equation}\label{moisa}
\beta(\Theta) \ge
\frac{\alpha(\Theta)}{2} \left( \frac{\alpha (\Theta)}{1-\alpha(\Theta)}
+\sqrt{\left(\frac{\alpha(\Theta)}{1-\alpha(\Theta)}\right)^2
+\frac{4\alpha(\Theta)}{1-\alpha(\Theta)}}\right) 
.\end{equation}
}

In \cite{m1}, \cite{m2} N. Moshchevitin obtained
the bounds in the cases $m=3, n=1$
and $m=n=2$. We will refer to a result from \cite{m2} (Theorem 24) which is the best know up to now in the case $m=3, n=1$.

{\bf Theorem E.} (N. Moshchevitin \cite{m1,m2})\,\,{\it
Suppose that $m=3,n=1$
and the collection $   \theta_1^1, \theta_1^2 , \theta_1^3
$ consists of numbers which, together with 1, are linearly
independent over $\mathbb{Z}$. Then
$$
\beta(\Theta) \ge\alpha(\Theta)\cdot\left(
\sqrt{\alpha (\Theta)+ \frac{1}{\alpha(\Theta)}-\frac{7}{4}}+\frac{1}{\alpha(\Theta)} - \frac{1}{2}\right).
$$}

To finish this section we would like to formulate a  result by
V. Jarn\'{\i}k from \cite{Cze} from which he deduces the inequality (\ref{dva}) of Theorem A.

{\bf Theorem F.} (V. Jarn\'{\i}k \cite{Cze})\,\,{\it
Suppose that $ n \ge 2$ , and that matrix
\begin{equation}\label{matt}
\Theta
=\left(
\begin{array}{cc}
\theta_1^1&\theta_1^2\cr
\vdots &\vdots \cr
\theta_n^1&\theta_n^2
\end{array}
\right)
\end{equation}
satisfy the condition (\ref{poss}).
Suppose that
a positive function $\psi (t) $ is such that
$$
\lim_{t\to +\infty} t \psi (t) = 0.
$$
Suppose that
$$
\psi^{\rm sup}_\Theta (t) \le \psi (t)
$$
for all $t$ large enough. 
Then there exist arbitrary large values of $t$ such that
$$
\psi^{\rm sup}_\Theta (t) \le\psi\left(\frac{1}{6t\psi(t)}\right).
$$}

Some related results are discussed in our recent surveys \cite{m2,mar}
and in the papers by M. Waldschmidt \cite{W} and O. German \cite{gAA,gMJCNT}.

{\bf 2. The result.}

We give few coments on the part (ii) of Jarn\'{\i}k's Theorem A. First of all we note that  the inequality 
({\ref{o2})
is better than the trivilal bound (\ref{ab}) in the case $\alpha (\Theta ) >2$ only.
However   $\alpha (\Theta) $ can attain any value from the interval
$\left[\frac{2}{n},+\infty \right]$.
So Theorem A gives nothing for the values of $\alpha (\Theta)$
in the interval  $\left[\frac{2}{n},2\right ]$.

As it was mentioned in the previous section,
it is possible to improve the inequality (\ref{dva}) in the case $n\ge 2$.
A proof of a  certain inequality better than (\ref{o2}) was sketched in \cite{m1,m2}  (Theorem 22 from \cite{m2}).
However the inequality from \cite{m1,m2} is very weak.
Moreover it is better than (\ref{dva}) in the range
 $ 1<\alpha (\Theta) < \left( \frac{1+\sqrt{5}}{2}\right)^2$ only.

In the present paper we get an inequality which improves the inequality (\ref{dva})
of Theorem A for all values of $\alpha (\Theta) >1$.
This inequality is better than that from \cite{m1,m2}.

Put
$$
G(\alpha )=
\begin{cases}
\frac{1}{2}\left(\alpha^2-\alpha+1+\sqrt{
(\alpha^2-\alpha+1)^2 +4\alpha^2(\alpha-1)}\right)\,\,\,\,\text{if}\,\,\,\,
1\le \alpha \le 2\cr
\frac{1}{2}\left(
\alpha^2-1+\sqrt{(\alpha^2-1)^2+4\alpha (\alpha-1)}\right)
\,\,\,\,\text{if}\,\,\,\,
\alpha \ge 2
\end{cases}
$$
and define
$$
g(\alpha ) = \frac{G(\alpha)}{\alpha}.
$$
Note that for $\alpha>1$ the value $g(\alpha)$ is the largest solution of the equation
\begin{equation}\label{eq}
\alpha g =\max (\alpha -1,1)+\frac{\alpha(\alpha-1)}{g - \alpha +1}.
\end{equation}
One can see that $g(1)=1$ and
$$
g (\alpha ) >\max (\alpha -1,1)
$$
for $ \alpha >1$.

Now we formulate the main result of the present paper.

{\bf Theorem 1.}\,\,{\it
Suppose  that $m=2$ and $ n\ge 3$.
Suppose that  among  $n+2$ two-dimensional vectors
\begin{equation}\label{vevel}
\left(
\begin{array}{c}
\theta_1^1\cr
\theta^2_1
\end{array}
\right),
\cdots ,
\left(
\begin{array}{c}
\theta_n^1\cr
\theta^2_n
\end{array}
\right),
\,\,\,
\left(
\begin{array}{c}
1\cr
0
\end{array}
\right),
\,\,\,
\left(
\begin{array}{c}
0\cr
1
\end{array}
\right)
\end{equation}
 there exist at least four vectors linearly independent over $\mathbb{Z}$.
Suppose that $\alpha (\Theta) \ge 1$.

Then
\begin{equation}\label{mai}
\beta(\Theta) \ge 
G(\alpha(\Theta)) =
\alpha(\Theta)  \cdot g(\alpha(\Theta)).
\end{equation}
}

{\bf Remark 1.}
From the conditions of Theorem it follows that for the matrix $\Theta$
one has (\ref{poss}).

{\bf Remark 2.}
The condition concerning linearly independence  of vectors (\ref{vevel}) cannot be removed in Theorem 1.
For example in the case $m=n=2$ one may 
take arbitrary $\alpha, \beta$ under the conditions
$\beta \ge \alpha(\alpha -1), \,\alpha \ge 2$ and
consider a matrix 
$$
\Theta=
\left(
\begin{array}{cc}
\theta_1^1& \theta_1^2\cr
\theta_1^1&\theta_1^2
\end{array}
\right)
$$
where $\theta_1^1,\theta_1^2$
come from Theorem B (ii).
Then $\alpha (\Theta) =\alpha,\,\,\beta(\Theta) =\beta$
and so (\ref{mai})
may be not true.

{\bf Remark 3.} Theorem 1 gives a bound which is better than the trivial bound (\ref{ab})  in the case $ \alpha (\Theta) > 1$ only.
In the case $ n=2$ we know that $\alpha (\Theta)$ cannot be less than one (see (\ref{miiha})).
However in the case $ n \ge 3$ we do not know if the trivial bound (\ref{ab}) can be improved upon in the range $
\frac{2}{n} < \alpha (\Theta ) < 1$.

{\bf Remark 4.}
If $ n \ge 3$ the trivial bound (\ref{ab}) cannot be 
 improved in the case $\alpha (\Theta ) = 1$, in general.
We refer to a result from 
 \cite{m2}
 (Theorem 10  and Corollary to it from \cite{m2}).
Suppose that $\xi \in \mathbb{R}\setminus \mathbb{Q}$ has bounded partial quotients in its continued fraction expansion.
Consider the matrix
\begin{equation}\label{matts}
\Theta
=\left(
\begin{array}{cc}
\theta_1&\xi \theta_1\cr
\theta_2 &\xi\theta_2\cr
\vdots &\vdots \cr
\theta_n&\xi\theta_n
\end{array}
\right),\,\,\,\, n \ge 3.
\end{equation}
 Then for almost all (in the sense og Lebesgue measure) real vectors $(\theta_1,\theta_1,...,\theta_n)\in \mathbb{R}^n$
all but a finite number of the best apppoximations vectors ${\bf Z}_\nu$ (see Sections 3,4 below) lie in a certain two-dimensional linear subspace of $\mathbb{R}^{n+2}$
and so for the matrix (\ref{matts}) one has $\alpha (\Theta) = 2$.
As the partial quotients of $\xi$ are bounded? one can see that $\beta (\Theta) = 1$ also.
Of course in this example all the elements of the matrix (\ref[matts}) can be linearly independent over $\mathbb{Z}$ together with 1.
So this example shows that for $n\ge 3$ it may happen that
$$
\beta (\Theta) = \alpha(\theta ) = 1,
$$
and the trivial bound (\ref{ab}) cannot be improved upon under the general condition of linear independence.
 
{\bf Remark 5.} In some very special cases (see the first Remark in  Section  6 below) it is possible to improve upon the trival buond (\ref{ab}) in the case $\alpha (\Theta) <1$.

{\bf 3. Ordinary best approximations.}

Recall the definition and the simplest properties of ordinary best approximation vectors.  These best approximations were actually used in the original paper  \cite{Cze} as well as in authors papers \cite{m1,m2,mCze,mar}.
As usual the  $\sup$-norm was used there to define the sequence of the best approximation vectors.

For  an integer vector  $ {\bf  x} = (x_1,x_2)\in \mathbb{Z}^2 $, put
$$
  \zeta^{\rm sup} ({\bf x}) =  \max_{1\le j\le n}||\theta_j^1x_1+\theta_j^2 x_2||.
$$
A vector
  $  {\bf x} \in \mathbb{Z}^2 $ is said to be  a {\it best approximation vector} if
$$
\zeta^{\rm sup} ({ \bf x})=\min_{{\bf  x}'} \zeta^{\rm sup} ({\bf x}'),$$ where the minimum is
taken over all  $ { \bf x}'  = (x_1', x_2')\in \mathbb{Z}^2 $ such that
$
0<| {\bf x}_i'|^{\rm sup}\le  | {\bf x}_i|^{\rm sup}. $

Suppose that the matrix $\Theta$ of the form (\ref{matt}) satisfies the following condition ({\bf L.I.}) : for any pair $(i,j), \,\, 1\le i ,j \le n,\,\, i \neq j$ the collection
$$
\theta_i^1,\theta_i^2,\theta_j^1,\theta_j^2, 1
$$
consists of numbers linearly independent over $\mathbb{Z}$.
From this condition on the matrix $\Theta$ 
we see that   all best approximations form the sequence
$${\bf x}^{\rm sup}_\nu =(x^{\rm sup}_{\nu,1},x^{\rm sup}_{\nu,2}),
\,\,\,\, \nu =1,2,3,...\,\, , $$
in such a way that for the values
$\zeta^{\rm sup}_\nu = \zeta^{\rm sup} ({\bf x}_\nu)$ and $X^{\rm sup}_\nu = |{\bf x}_{\nu,i}|^{\rm sup}$
form infinite monotone sequences
\begin{equation}\label{11}
\zeta^{\rm sup}_1> \zeta^{\rm sup}_2>...>\zeta^{\rm sup}_\nu>\zeta^{\rm sup}_{\nu+1}>... \,\,\, , \end{equation}
\begin{equation}\label{22}
X^{\rm sup}_1< X^{\rm sup}_2<...<X^{\rm sup}_nu<X^{\rm sup}_{\nu+1}<... \,\,\, .
\end{equation}
For a best approximation vector $ {\bf x}^{\rm sup}_\nu = (x^{\rm sup}_{\nu,1}x^{\rm sup}_{\nu,2})$
we consider integers $y^{\rm sup}_{\nu,j},\, 1\le j\le n$ defined by the equalities
$$
||\theta_j^1x^{\rm sup}_{\nu,1}+\theta_j^2
x^{\rm sup}_{\nu,2}||=
|\theta_j^1x^{\rm sup}_{\nu,1}+\theta_j^2
x^{\rm sup}_{\nu,2}- y^{\rm sup}_{\nu,j}|
$$
and define the {\it extended best approximation vector}
$$
{\bf z}^{\rm sup}_\nu 
=(x^{\rm sup}_{\nu,1},x^{\rm sup}_{\nu,2}, y^{\rm sup}_{\nu,1},....,y^{\rm sup}_{\nu,n}) \in \mathbb{Z}^{n+2}.
$$
Here we should note that each vector ${\bf z}^{\rm sup}_\nu$ is a primitive vector, that is
$$
{\rm g.c.d.}
(x^{\rm sup}_{\nu,1},x^{\rm sup}_{\nu,2}, y^{\rm sup}_{\nu,1},....,y^{\rm sup}_{\nu,n}) =1.
$$
Moreover each couple of consecutive vectors ${\bf z}^{\rm sup}_\nu, {\bf z}^{\rm sup}_{\nu+1}$
can be extended to a basis of the whole integer lattice $\mathbb{Z}^{n+2}$.
In particular ${\bf z}^{\rm sup}_\nu$ and ${\bf z}^{\rm sup}_{\nu+1}$
are linearly independent.

We consider two-dimensional subspace
$$
{\cal L} =\{
(x_1,x_2,y_1,...,y_n) \in \mathbb{R}^{n=2}:\,\,
\theta_j^1x_1+\theta_j^2 x_2 -y_j =0,\,\, 1\le j \le n\}.
$$
From ({\bf L.I.}) condition   on the matrix $\Theta$  
  it follows that there is no non-zero integer points in ${\cal L}$ and the best approximation vectors ${\bf z}^{\rm sup}_\nu$
become more and more close to ${\cal L}$ as $\nu$ tends to infinity.
From Minkowski convex body theorem it follows that
\begin{equation}\label{mii}
\zeta^{\rm sup}_\nu (X^{\rm sup}_{\nu+1})^{\frac{2}{n}}\le 1.
\end{equation}
Here we should note that the inequality
$$
\psi^{\rm sup}_\Theta (t) \le t^{-\alpha}
$$
holds for al $t$ large enough if and only if
$$
\zeta^{\rm sup}_\nu  \le (X^{\rm sup}_{\nu+1})^{-\alpha} $$
for $\nu$ large enough.

{\bf 4. Spherical best approximations.}

However consideration of the ordinary best approximations vectors is not very convenient for our purposes. It makes the proofs too
cumbersome.
To make
our proofs easier we need another definition. 

In the sequel by ${\rm dist } ({\cal A}, {\cal B})$ we denote the Euclidean distance between the sets ${\cal A}, {\cal B}\subset \mathbb{R}^{n+2}$.
We shall consider vectors from $\mathbb{R}^{n+2}$ of the form
$$
{\bf z} = (x_1,x_2, y_1,...,y_n).
$$
For such a vector by $Z = Z({\bf z}) =  {\rm dist} (\{ {\bf z}\}, \{{\bf 0}\})$
 we define its Euclidean norm and by
$\zeta ({\bf z})
={\rm dist}( \{{\bf z}\}, {\cal L})
$
 we define the distance from ${\bf z}$ to 
the two dimensional subspace ${\cal L}$
defined in the previous section.

We need a simple geometric observation.

{\bf Lemma 1.}
\,\,{\it
Let $\pi$ be a two-dimrnsional linear subspace in $\mathbb{R}^{n+2}$ such that
$\pi \cap {\cal L} = \{ {\bf 0}\}$. Suppose that $\pi$ and ${\cal L}$ are not orthogonal.
Then given $\lambda >0$ the set
$${\cal G}_\lambda =
\{ {\bf z} \in \pi:\,\, 
{\rm dist } (\{{\bf z}\}, {\cal L}) = \lambda
	\}
$$
is an ellipse.
Moreover for
all values of $\lambda$  all the ellipses ${\cal G}_\lambda$
are dilatated form the 
ellipse ${\cal G}_1$, and hence all their minor axes coinside and all their major
axes coinside.}

{\bf Remark.}\,\,
It is clear that in the case ${\rm dim }\, \pi \cap {\cal L}  = 1$
the set ${\cal G}_\lambda$ consists of two parallel lines.

Proof of Lemma 1.

We may restict ourselves on four-dimensional  subspace
${\rm span}\, (\pi\cup{\cal L})$.

Suppose that $\eta_1, \eta_2$ are linearly independent vectors from ${\cal L}$.

Then the Euclidean distance from ${\bf z}\in \pi $ to ${\cal L}$
is defined by the formula
$$
{\rm dist }\, (\{{\bf z}\}, {\cal L}) =
\frac{\text{ volume of the parallelepiped spaned by}\,\, {\bf z}, \eta_1,\eta_2}{
\text{area of the parallelogram spaned by}\,\, \eta_1,\eta_2}.
$$
So ${\rm dist }\, (\{{\bf z}\}, {\cal L})$ is a quadratic form in ${\bf z}$.
Being restricted on $\pi$ it gives a quadratic form in two variables.
It is clear that ${\cal G}_\lambda$ is a bounded set. So ${\cal G}_\lambda$ is an ellipse. 
Further statements of Lemma 1 are 
obvious.$\Box$

For a two-dimensional linear subspece $\pi $ the following observation will be of importance.
Consider the circle 
$$
\hbox{\got S}=\{ {\bf z}\in \pi :\,\,\,
Z({\bf z}) = 1\}.
$$
Let $\pi$ be not an orthogonal complement to ${\cal L}$.
Then there exist two orthogonal vectors ${\bf a}, {\bf b}\in \pi $
such that 
$$
\min_{{\bf z}\in \hbox{\got S}}
\,{\rm dist} \,(\{{\bf z}\}, {\cal L}) =
{\rm dist}\, (\{{\bf p}\}, {\cal L})
,\,\,\,
\max_{{\bf z}\in \hbox{\got S}}
\,{\rm dist} \,(\{{\bf z}\}, {\cal L}) =
{\rm dist}\, (\{{\bf q}\}, {\cal L})
.
$$
We supppose that the directed angle between vectors ${\bf p}$ and ${\bf q}$ is equal to $+\pi/2$.
For $t\in [0.\pi/2]$  we  consider the point ${\bf  A}(t)\in \hbox{\got S}$ obtained by the rotation of the point ${\bf p}\in \hbox{\got S}$
by the angle $t$ towards the point ${\bf q}\in \hbox{\got S}$.
We are interested in the function
$$
f(t) =\frac{ {\rm dist} \,(\{{\bf A}(t)\}, {\cal L})}{{\rm dist} \,(\{{\bf A}(t)\}, {\rm span} \, {\bf p)}} =
\frac{ {\rm dist} \,(\{{\bf A}(t)\}, {\cal L})}{\sin t}
.
$$

{\bf Lemma 2.}\,\,{\it
In the interval $0< t\le \pi/2$ the function $f(t)$ decreases.}

Proof.

We may suppose that  both two-dimensional subspaces  belong to the same four-dimensional Euclidean subspace
$\mathbb{R}^4$ with coordinares $\eta_1,\eta_2,\eta_3,\eta_4$ and that the subspace ${\cal L}$
in these coordinates is determined by the equations
$$
\eta_1=\eta_2 = 0.
$$
If in these coordinates we have a point $ {\bf z} = (\eta_1,\eta_2,\eta_3,\eta_4)$,
then
$$
{\rm dist} \,(\{{\bf z}\}, {\cal L}) =\sqrt{\eta_1^2+\eta_1^2}.
$$
Let in these coordinates 
$$
{\bf p} = (p_1,p_2,p_3,p_4),\,\,\,\,
{\bf q} = (q_1,q_2,q_3,q_4).
$$
Then
$$
(f(t))^2 = 
\frac{
(p_1\cos t+q_1\sin t)^2 +(p_2\cos t +q_2\sin t)^2
}{\sin^2 t} =
 \frac{p_1^2+p_2^2}{\sin^2 t} + q_1^2+q_2^2 -p_1^2-p_2^2
$$
(we should note that the point ${\bf p}$ is the closest point to ${\cal L}$ and so
$ p_1q_1+p_2 q_2 =0$).
Now it is clear that $f(t)$ decreases.$\Box$

 {\bf Remark.}\,\, In the case
${\rm dim}\, \pi \cap {\cal L} = 1$
the function $f(t)$ is a constant as in this case $ p_1^2+p_2^2 =0$.

Now we define the sequence of {\it spherical best approximation vectors}.

We define ${\bf z}\in \mathbb{Z}^{n+2}$ to be a spherical best approximation vector
if
$$
\zeta ({\bf z}) \le \zeta ({\bf z}')
$$
for all nonzero integer vectors ${\bf z}'\in \mathbb{Z}^{n+2}$
with $ Z' =  {\rm dist} \,( {\bf z}', \{{\bf 0}\}) \le Z$.

To avoid the situation when two best approximation vectors with the same value of $Z$ may occur we need to suppose a condition which generalizes the condition
({\bf L.I.}) from the previous section.
However such a condition deal with quadratic relations instead of linear relations.  We do not want to suppose additional restrictions on matrix $\Theta$.
So we will not define the sequence of the best spherical approximation vectors in a unique way. 
Analoguosely to the sequences of the ordinary best approximations ${\bf z}^{\rm sup}_\nu$ satisfying
(\ref{11},\ref{22}) we define the sequence
\begin{equation}\label{000}
	{\bf z}_\nu = (x_{\nu,1}, x_{\nu,2}, y_{\nu,1},...,y_{\nu,n}),\,\,\,\,
\nu =1,2,3,...
\end{equation}
such that for 
$\zeta_\nu = \zeta ({\bf z}_\nu)$ and 
$Z_\nu = Z({\bf z}_\nu)$
one has
\begin{equation}\label{111}
\zeta_1> \zeta_2>...>\zeta_\nu>\zeta_{\nu+1}>... \,\,\, , \end{equation}
\begin{equation}\label{222}
Z_1< Z_2<...<Z_\nu<Z_{\nu+1}<... \,\,\, .
\end{equation}
Of course under the conditions of Theorem 1 
it may happen that the same values of $\zeta_\nu, Z_\nu$
are attained on two (or even more) different integer vectors.
In such a situation we choose one of the admisssible integer vectors
in an arbitrary way and define it to be the $\nu$-th best spherical approximation
vector.
So the sequence (\ref{000})
may depend on our choice.  But the values from the sequences (\ref{111},\ref{222}) do not depend on our choice.
Aftrer we have chosen the sequence (\ref{000}) we fix it.
Everywhere in the sequel we deal with the  fixed seuqence of spherical best approximations which was chosen here.

Analogously to the ordinary best approximations,
 each vector ${\bf z}_\nu$ is a primitive vector,  and
each couple of consecutive vectors ${\bf z}_\nu, {\bf z}_{\nu+1}$
can be extended to a basis of the whole integer lattice $\mathbb{Z}^{n+2}$,
and in particular ${\bf z}_\nu$ and ${\bf z}_{\nu+1}$
are linearly independent.

The set
$$
\{ {\bf z} \in \mathbb{R}^{n+2}:\,\,\,
\zeta ({\bf z}) <\zeta_\nu,\,\, Z({\bf z}) < Z_{\nu+1}\}
$$
has no non-zero integer points inside.
So analogously to (\ref{mii})
from Minkowski convex body theorem we have
\begin{equation}\label{mii1}
\zeta_\nu Z_{\nu+1}^{\frac{2}{n}}\ll_n 1,
\end{equation}
where the constant in the symbol $\ll_n$ may depend on the dimension $n$.

Put
$$
\psi _\Theta (t) =
 \min_{{\bf z} \in \mathbb{Z}^{n+2}:\,0< Z \le t}\,\,\, 
{\rm dist} \, (\{{\bf z}\}, {\cal L}).
$$
As all the norms in Euclidean spaces are equivalent, we see that
in the definitions (\ref{aa},\ref{bb}) of the exponents $\alpha (\Theta)$ and $\beta(\Theta)$
we may replace
the function 
$\psi_\Theta^{\rm sup} (t)$ by the function
$\psi_\Theta (t)$ and the result will be the same. So
$$
\psi_\Theta (t) \le t^{-\alpha}$$
for all $t$ lagre enough if and only if
$$
\zeta_\nu \le Z_{\nu+1}^{-\alpha}
.$$
In particular for any $\alpha <\alpha(\Theta)$ for all $\nu$ large enough one has
\begin{equation}\label{useless}
\zeta_\nu < Z_{\nu+1}^{-\alpha}.
\end{equation}

{\bf 5.  Successive best approximation vectors in two-dimensional subspace.}

It may happen that three or more vectors 
\begin{equation}\label{nuk}
{\bf z}_\nu, {\bf z}_{\nu+1},...,{\bf z}_k
\end{equation}
lie in a cetrain two-dimensional linear subspace $\pi \subset \mathbb{R}^{n+2}$.
Then the following statement is valid.

{\bf Lemma 3.}\,\,{\it In the case when vectors (\ref{nuk}) lie in a certain  two-dimensional linear subspace $\pi$ one has
\begin{equation}\label{norma}
\zeta_\nu Z_{\nu+1} \le \frac{12}{\sqrt{\sqrt{52}-5}}\cdot \zeta_{k-1}Z_k.
\end{equation}}

{\bf Corollary.}\,\,{\it
Suppose that $\alpha <\alpha (\Theta)$.
Then if $\nu$ is large enough and
 vectors (\ref{nuk}) lie in a certain  two-dimensional linear subspace $\pi$ one has
\begin{equation}\label{norma1}
\zeta_\nu \ll Z_{\nu+1}^{-1}Z_{k}^{1-\alpha}.
\end{equation}}

Proof of Lemma 3.

First of all we consider the case $ \pi \cap {\cal L} = \{ {\bf 0}\}.$

Consider the collection of ellipses  $\{ {\cal G}_\lambda\}_{\lambda>0}
$
They have  commom major axes. We denote the one-dimensional subspace of major axes by ${\cal P}$.
The orthogonal one-dimensional subspace 
consisting of all common
minor axes
we denote by ${\cal Q}$.

For every $l$ from the interval $\nu\le l\le k-1$ we define the value of $\lambda_l$
from the condition
$$
{\bf z}_l \in {\cal G}_{\lambda_l}.
$$
As
$$
\zeta_\nu >\zeta_{\nu+1}>...>\zeta_{k-1}>\zeta_k
$$
we have
$$
\lambda_\nu >\lambda_{\nu+1}>...>\lambda_{k-1}>\lambda_k.
$$
For $\nu\le l\le k-1$ we put
$$
\xi_l =  {\rm dist}\, (\{{\bf z}_l\}, {\cal P}).
$$
One can see  by the monotonicity argument  that
$$
\xi_\nu >\xi_{\nu+1}>...>\xi_{k-1}>\xi_k.
$$ 

Define $\Xi_l$ to be the length of a half of the minor axis of the ellipse
${\cal G}_{\lambda_l}$.
Then
\begin{equation}\label{elle}
\Xi_l\ge \xi_l
.
\end{equation}
The planar convex set
$$
{\cal E}_l=
\{ {\bf z} \in \pi:\,\,\,
\zeta ({\bf z}) <\zeta_l,\,\, Z({\bf z}) < Z_{l+1}\}
\subset \pi$$
has no non-zero integer points inside.
There are two pairs of independent  integer points $\pm {\bf z}_l, \pm{\bf z}_{l+1}$
on its boundary.
For the two-dimensional volume of the set ${\cal E}_l$ one has  estimates
$$
2\Xi_lZ_{l+1} \le {\rm vol}_2 \,{\cal E}_l
\le 4\Xi_lZ_{l+1} .
$$
Consider the two-dimensional lattice
$$\Lambda =
\mathbb{Z}^{n+2}\cap \pi
$$
with the two-dimensional fundamental volume ${\rm det}\,\Lambda$.
Then by Minkowski convex  body theorem
$$
\frac{\Xi_lZ_{l+1}}{2}\le {\rm det}\, \Lambda
\le 2\Xi_lZ_{l+1}.$$
So
\begin{equation}\label{vva}
\Xi_\nu Z_{\nu+1} \le 4 \Xi_{k-1}Z_k.
\end{equation}

Put
$$ a=
\frac{\sqrt{\sqrt{52}-5}}{3}
,
\,\,\,\,
b =\sqrt{1-a^2}
.
$$
Now we   prove the inequalities
\begin{equation}\label{monot}
\xi_l\ge a\Xi_l,\,\,\,\,
\nu\le l \le k-1.
\end{equation}
Indeed there are at least two independent points ${\bf z}_{l-1}, {\bf z}_l \in \Lambda$ on the boundary of the set ${\cal E}_{l-1}$ . Thus 
${\bf z}_{l-1}
$
lies on the boubdary of the ellipse
${\cal G}_{\lambda_{l-1}}$
and
${\bf z}_l$
lie inside the ellipse ${\cal G}_{\lambda_{l-1}}$.
 
Let $H$ be the half of the lenght of the major axis of  ${\cal G}_{\lambda_{l-1}}$.
If the distance from ${\bf z}_{l-1} $ to the minor axis of
${\cal G}_{\lambda_{l-1}}$ is greater than
$bH$ then the
the distance from ${\bf z}_{l} $ to the minor axis of
${\cal G}_{\lambda_{l-1}}$ is greater than
$bH$ also. Easy calculation shows that in this case the  point
$ {\bf z}_{l-1} -{\bf z}_l$ lies inside ${\cal E}_{l-1}$.
It is not possibe.
So the distanse from
 ${\bf z}_{l-1} $ to the minor axis of
${\cal G}_{\lambda_{l-1}}$ is not greater  than
$bH$.
Hence the distance from
 ${\bf z}_{l-1} $ to the major axis of
${\cal G}_{\lambda_{l-1}}$ is not less than
$a\Xi_{l-1}$. But this distance is equal to $\xi_{l-1}$.
So we have
$$
\xi_{l-1} \ge a \Xi_{l-1}.
$$
Inequalities (\ref{monot}) are proved.

Applying the ${\bf 0}$-central projection onto $\hbox{\got S}$
and taking into account Lemma 2
we see that
\begin{equation}\label{corol}
\frac{\zeta_{k-1}}{\xi_{k-1}}\ge \frac{\zeta_\nu}{\xi_\nu}.
\end{equation}

From (\ref{elle},\ref{vva},\ref{monot},\ref{corol}) we immediately deduce (\ref{norma}).

The case $ {\rm dim}\, \pi \cap {\cal L} = 1$ is easier\footnote{In the sketched proof  in \cite{m1,m2} only this case was considered.}.
In this case the set ${\cal G}_\lambda$ is a union of two parallel lines and the function $f(t)$ is a constant function.
So
we may  assume that $\xi_l = \Xi_l, \nu\le l\le k-1$,
inequalities (\ref{vva}) remain true, and instead of (\ref{corol}) one has
$\frac{\zeta_{k-1}}{\xi_{k-1}}= \frac{\zeta_\nu}{\xi_\nu}$.
So  (\ref{norma}) follows in this case also.$\Box$

{\bf Remark to the proof of Lemma 3.}\,\,
Similar argument was used not by the author in \cite{m1,m2,mCze} only,
but by some other mathematicians.
In particular similar argument was applied by Y. Cheung  (see Theorem 1.6 from  \cite{cheu}).

{\bf 6. Dimension of subspace of  best approximation vectors.}

 It may happen that all the best approximation vectors
$ {\bf z}_\nu$ lie in a certain linear subspace of $\mathbb{R}^{n+2}$
of dimension less than $n+2$. So we consider the value
$$
R(\Theta) =
\min \{ r:\,\,
\text{there exist a subspace}\,\, {\cal R}\subset\mathbb{R}^{n+2}
\,\,\text{and}\,\, \nu_0\in \mathbb{Z}_+\,\,
\text{s.t.}\,\, \forall\, \nu \ge \nu_j\,\,\, {\bf z}_\nu \in {\cal R}\}.
$$ 
Let ${\cal R} ={\cal R}(\Theta)$ be the linear subspace from the definition of $R(\Theta)$. We see that the lattice $\mathbb{Z}^{n+2}\cap {\cal R}$ 
is a lattice of dimension $R(\Theta)$.
Here we would like to recall the definition of {\it completely rational } subspace.
A subspace $\pi \subset \mathbb{R}^d$ is defined to be
{\it completely rational} if the lattice $\mathbb{Z}^d\cap \pi$ has dimension
$d$. So ${\cal R}$ is a completely rational subspace.

Consider the subspace
$$
{\cal K} = {\cal K}(\Theta) ={\cal R}\cap {\cal L}.
$$
As there is an infinite sequence of integer points ${\bf z}_\nu \in {\cal R}$
such that the distance between ${\bf z}_\nu$ and ${\cal L}$ tends to zero as $\nu$ tends  to infinity,
we see that ${\rm dim } \, {\cal K} >0$.
So we have only two opportunities.
Either  
${\rm dim }\, {\cal K} = 1$, or
${\rm dim } \, {\cal K}=2$ and in this case $ {\cal K} = {\cal L} \subset {\cal R}$.

The case ${\rm dim}\, {\cal K} = 1$ is easy.
First of all we note that there is no completely rational subspace $\pi 
\subset \mathbb{R}^{n+2}$ such that
${\cal K}\subset \pi\subset {\cal R}$
and ${\rm dim}\, \pi < {\rm dim}\,{\cal R}$.

If ${\rm dim}\, {\cal R} = 2$ then we deal with approximations to a one-dimensional  subspace ${\cal K}$ from the two-dimensional rational subspace ${\cal R}$.
It this case $\alpha(\Theta) = 1$ and there is nothing to prove.

If ${\rm dim}\, {\cal R} > 2$ then
  we deal with the approximations to
a one-dimensional  subspace ${\cal K}$ from the  rational subspace ${\cal R}$
of dimension greater than two. 
In this case
$$
\frac{1}{{\rm dim}\, {\cal R} -1} \le \alpha(\Theta)\le 1.
$$
So this case does not considered in Theorem 1.

{\bf Remark.}\,\,
The situation in the case
${\rm dim}\, {\cal K} = 1$, ${\rm dim}\, {\cal R} > 2$
 is quite similar to the setting which was considered in Theorem A,
statement (i). So in this case the inequality (\ref{sim}) is valid. The inequality (\ref{sim})  gives an optimal boumd in the case $ {\rm dim}\, {\cal R}=3$.
  If 
${\rm dim}\, {\cal R}>3$ we may apply Theorem C part (i) or Theorem D and obtain an even better lower bound for $\beta (\Theta)$ in terms of $\alpha (\Theta)$:
for $R(\Theta) = 4$  we have the bound  (\ref{moisa}) from Theorem D, and  for $R(\Theta) > 4$ from 
the inequality (\ref{o1}) of Theorem C, part (i) we have
$$
\beta (\Theta)  \ge \alpha(\Theta ) \cdot
\frac{\alpha (\Theta)+ {\rm dim}\, {\cal R}-3 }{({\rm dim}\, {\cal R}-2)(1-\alpha (\Theta))}
.
$$
In any case here we have a bound which is better than the trivial bound (\ref{ab}).

Now we consider the case ${\rm dim } \, {\cal K}=2$  when 
 $ {\cal K} = {\cal L} \subset {\cal R}$.
Here we should not that under the conditions
of  Theorem 1 
among two-dimensional vectors (\ref{vevel}) there are at least four  vectors
  linearly independent over $\mathbb{Z}$.
So ${\cal L} $ cannot lie in a completely rational subspace of dimension $\le 3$,
and if ${\rm dim}\, {\cal K}(\Theta) = 2$ then
\begin{equation}\label{4}
{\rm dim }\, {\cal R}(\Theta) \ge 4.
\end{equation} 
(For more details see Section 2.1 from \cite{m2} and especially formula (21).)
In the rest of the paper we suppose that (\ref{4}) holds.
 
{\bf Remark.}\,\
The author does not know if in the case $m=n = 2$ there exists a
 matrix
$\Theta$ satisfying the condition of Theorem 
1  and such that ${\rm dim}\, {\cal R} (\Theta) = 2, \,\,{\rm 
dim }\, {\cal K} (\Theta) = 1.$
From Jarn\'{\i}k 's result it follows that in the case $m=n=2$
the situation  with ${\rm dim}\, {\cal R} (\Theta) = 3, \,\,{\rm 
dim }\, {\cal K} (\Theta) = 1$ never happens (see the discussion in \cite{m2,mar}).

{\bf 7. Four linearly independent vectors.}

 From the condition (\ref{4}) we see that
   there exist infinitely many pairs of indices
   $\nu<k, \nu\to +\infty$   such that

{\bf (a)}  both triples
$$
{\bf z}_{\nu-1},{\bf z}_\nu,{\bf z}_{\nu+1};\,\,\,\,\,\, {\bf
z}_{k-1},{\bf z}_k,{\bf z}_{k+1}
$$ consist of
  linearly independent vectors;

{\bf (b)} there exists a two-dimensional linear subspace   $\pi$
such that
$$
{\bf z}_l\in \pi,\,\,\, \nu\le l\le k;\,\,\,\,\, {\bf z}_{\nu-1}
\not\in \pi,\,\,\, {\bf z}_{k+1} \not\in \pi;
$$

{\bf (c)}  the vectors
$$
{\bf z}_{\nu-1},{\bf z}_\nu,{\bf z}_{k},{\bf z}_{k+1}
$$
are linearly independent.

If a pair if indices $(\nu,k)$ satisfy  {\bf (a)}, {\bf (b)}, {\bf(c)} we 
say that $(\nu,k)$ satisfy  {\bf (abc)}-property.

Lemmas 4 and 6   below   were actually proved by Jarn\'{\i}k in \cite{Cze}.
However they were not stated by him explicitly. So we give a complete proof.
Lemma  5 comes  from \cite{m1,m2}.

{\bf Lemma 4.}\,\,
{\it
Suppose that $\alpha (\Theta)>1$. Then for all $\nu$ large enough one has
$$
\Delta_\nu =
\left|
\begin{array}{cc}
x_{\nu,1} & x_{\nu,2}\cr
x_{\nu+1,1}&x_{\nu+1,2}
\end{array}
\right|
\neq 0.
$$
}

Proof.
Suppose that
$\Delta_\nu = 0$.
Consider the determinants
$$
\Delta_{\nu,j } =
\left|
\begin{array}{cc}
x_{\nu,1} & y_{\nu,j}\cr
x_{\nu+1,1}&y_{\nu+1,j}
\end{array}
\right|
 =
\left|
\begin{array}{cc}
x_{\nu,1} 
& y_{\nu,j} - \theta^1_jx_{\nu,1} -\theta^2_j x_{\nu,2}
\cr
x_{\nu+1,1}
&
y_{\nu+1,j} - 
\theta^1_jx_{\nu+1,1}-\theta^2_jx_{\nu+1,2}
\end{array}
\right|
,\,\,\,\,
1\le j \le n.
$$
 As $\alpha (\Theta) >1$ we see that
$$
|\Delta_{\nu,j} |
\le 2 Z_{\nu+1}\zeta_\nu \to 0,\,\,\,\, \nu \to \infty.
$$
That is why
$$
\Delta_{\nu,j} = 0,\,\,\,\, 1\le j \le n.
$$
But we have supposed that $\Delta_\nu = 0$ also.
This means that the vectors ${\bf z}_\nu$ and ${\bf z}_{\nu+1}$ are linearly dependent. This is a contradiction. $\Box$

{\bf Lemma 5.}\,\,{\it Suppose that $\alpha (\Theta)>1$ and $1<\alpha <\alpha(\Theta)$.
Suppose that the pair of indices $(\nu,k)$ satisfies {\bf (abc)}-property and $\nu$ is large enough. Then
 \begin{equation}\label{ge1}
Z_{k+1}\gg Z_\nu^\alpha Z_k^{\alpha - 1}.
\end{equation}
}

Proof.

To prove (\ref{ge1}) we consider four linearly independent integer vectors
${\bf z}_{\nu-1},{\bf z}_{k-1}, {\bf z}_{z}, {\bf z}_{k+1}$.
Consider four-dimensional vectors
$$
\left(
\begin{array}{c}
x_{\nu-1,1}\cr
x_{k-1,1}\cr
x_{k,1}\cr
x_{k+1,1}
\end{array}\right)
,
\left(
\begin{array}{c}
x_{\nu-2,1}\cr
x_{k-1,2}\cr
x_{k,2}\cr
x_{k+1,2}
\end{array}\right),
\left(
\begin{array}{c}
y_{\nu-1,1}\cr
y_{k-1,1}\cr
y_{k,1}\cr
y_{k+1,1}
\end{array}\right)
,
\cdots,
\left(
\begin{array}{c}
y_{\nu-1,n}\cr
y_{k-1,n}\cr
y_{k,n}\cr
y_{k+1,n}
\end{array}\right)
.
$$
Among these vectors there are four linearly independent ones.
From Lemma 4 we know that two vectors
$$
\left(
\begin{array}{c}
x_{\nu-1,1}\cr
x_{k-1,1}\cr
x_{k,1}\cr
x_{k+1,1}
\end{array}\right)
,
\left(
\begin{array}{c}
x_{\nu-2,1}\cr
x_{k-1,2}\cr
x_{k,2}\cr
x_{k+1,2}
\end{array}\right)
$$
are linearly independent.
So there exist indices $i\neq j$ such that four vectors
$$
\left(
\begin{array}{c}
x_{\nu-1,1}\cr
x_{k-1,1}\cr
x_{k,1}\cr
x_{k+1,1}
\end{array}\right)
,
\left(
\begin{array}{c}
x_{\nu-2,1}\cr
x_{k-1,2}\cr
x_{k,2}\cr
x_{k+1,2}
\end{array}\right),
\left(
\begin{array}{c}
y_{\nu-1,i}\cr
y_{k-1,i}\cr
y_{k,i}\cr
y_{k+1,i}
\end{array}\right)
,
\left(
\begin{array}{c}
y_{\nu-1,j}\cr
y_{k-1,j}\cr
y_{k,j}\cr
y_{k+1,j}
\end{array}\right)
.
$$
We consider the  determinant 
$$D
=
\left|
\begin{array}{cccc}
x_{\nu-1,1}&x_{\nu-1,2}&y_{\nu-1,i}&y_{\nu-1,j}\cr
x_{k-1,1}&x_{k-1,2}&y_{k-1,i}&y_{k-1,j}\cr
x_{k,1}&x_{k,2}&y_{k,i}&y_{k,j}\cr
x_{k+1,1}&x_{k+1,2}&y_{k+1,i}&y_{k+1,j}
\end{array}
\right|
=
$$
$$
=
\left|
\begin{array}{cccc}
x_{\nu-1,1}&x_{\nu-1,2}&y_{\nu-1,i} -   \theta^1_ix_{\nu-1,1}-\theta^2_ix_{\nu-1,2}&y_{\nu-1,j}  -   \theta^1_jx_{\nu-1,1}-\theta^2_jx_{\nu-1,2}\cr
x_{k-1,1}&x_{k-1,2}&y_{k-1,i} -   \theta^1_ix_{k-1,1}-\theta^2_ix_{k-1,2}&y_{k-1,j}  -   \theta^1_jx_{k-1,1}-\theta^2_jx_{k-1,2}\cr
x_{k,1}&x_{k,2}&y_{k,i} -   \theta^1_ix_{k,1}-\theta^2_ix_{k,2}&y_{k,j}  -   \theta^1_jx_{k,1}-\theta^2_jx_{k,2}\cr
x_{k+1,1}&x_{k+1,2}&y_{k+1,i} -   \theta^1_ix_{k+1,1}-\theta^2_ix_{k+1,2}&y_{k+1,j}  -   \theta^1_jx_{k+1,1}-\theta^2_jx_{k+1,2}\cr
\end{array}
\right|
$$
 corresponding to these  integer vectors.
As $D\neq 0$ we see that
$$
1\le |D| \le 24 \zeta_{\nu-1}\zeta_{k-1}Z_kZ_{k+1}\ll
Z_\nu^{-\alpha} Z_k^{1-\alpha} Z_{k+1}
$$
(here we use (\ref{useless})) and Lemma 5 follows.$\Box$

{\bf Lemma 6.}
\,\,
{\it 
Suppose that $\alpha (\Theta) >2. $
Then for any positive $\varepsilon$
for all $\nu$
 large enough if three best approximation vectors 
${\bf z}_{\nu-1},{\bf z}_\nu, {\bf z}_{\nu+1}$ are linearly independent then $Z_{\nu+1} \ge Z_\nu^{\alpha(\Theta ) -\varepsilon -1}$.}

Proof.

Consider three-dimensional vectors
\begin{equation}\label{tree}
\left(
\begin{array}{c}
 x_{\nu-1,1}\cr
x_{\nu,1} \cr
x_{\nu+1,1}
\end{array}
\right),
\,\,\,
\left(
\begin{array}{c}
 x_{\nu-1,2}\cr
x_{\nu,2} \cr
x_{\nu+1,2}
\end{array}
\right),
\,\,\,
\left(
\begin{array}{c}
 y_{\nu-1,1}\cr
y_{\nu,1} \cr
y_{\nu+1,1}
\end{array}
\right),
 \cdots ,
\left(
\begin{array}{c}
 y_{\nu-1,n}\cr
y_{\nu,n} \cr
y_{\nu+1,n}
\end{array}
\right).
\end{equation}
 As three vectors
 ${\bf z}_{\nu-1},{\bf z}_\nu, {\bf z}_{\nu+1}$ are linearly independent
we see that among three-dimensional vectors (\ref{tree}) there
are three linearly independent vectors.
 From Lemma 4  we know that three-dimensional vectors
$$
\left(
\begin{array}{c}
 x_{\nu-1,1}\cr
x_{\nu,1} \cr
x_{\nu+1,1}
\end{array}
\right)
\,\,\,
\text{and}
\,\,\,
\left(
\begin{array}{c}
 x_{\nu-1,2}\cr
x_{\nu,2} \cr
x_{\nu+1,2}
\end{array}
\right)
$$
are linearly independent.
So there exists $j$ such that
$$
\Delta =
\left|
\begin{array}{ccc}
x_{\nu-1,1} &x_{\nu-1,2}&y_{\nu-1,j}\cr
x_{\nu,1}&x_{\nu,2}&y_{\nu,j}\cr
x_{\nu+1,1} &x_{\nu+1,2}&y_{\nu+1,j}
\end{array}
\right|
=
\left|
\begin{array}{ccc}
x_{\nu-1,1} &x_{\nu-1,2}&y_{\nu-1,j}-\theta^1_jx_{\nu-1,1}-\theta^2_jx_{\nu-1,2}\cr
x_{\nu,1}&x_{\nu,2}&y_{\nu,j} -\theta^1_jx_{\nu,1}-\theta^2_jx_{\nu,2}\cr
x_{\nu+1,1} &x_{\nu+1,2}&y_{\nu+1,j}-\theta^1_jx_{\nu+1,1}-\theta^2_jx_{\nu+1,2}
\end{array}
\right|
\neq 0.
$$
Now we consider the inequality
$$
1\le |\Delta | \le 6Z_{\nu+1}Z_\nu \zeta_{\nu-1}.
$$
But for $\nu$ large enough we have $ \zeta_{\nu-1} \le Z_\nu^{-\alpha(\Theta) +\varepsilon}$,
and Lemma 6 follows.$\Box$

{\bf  8. Proof of Theorem 1.}

We take $\alpha<\alpha (\Theta)$ close to $\alpha (\Theta)$.
Suppose that $(\nu, k)$ satisfies {\bf (abc)}-property and $\nu$ is large enough.

If $\alpha (\Theta) >2$ and
$\alpha >2$ from Lemma 6
we have
$
Z_{\nu+1} \gg Z_{\nu}^{\alpha -1}.
$
If $\alpha (\Theta ) \le 2$ we have nothing but trivial bound
$
Z_{\nu+1} > Z_{\nu}.
$
So in any case
\begin{equation}\label{u2}
Z_{\nu+1} \gg Z_{\nu}^{\max(\alpha -1,1)}.
\end{equation}

Now we prove that either
$$
Z_{\nu+1} \gg Z_{\nu}^{g(\alpha)}
$$
or
$$
Z_{k+1} \gg Z_{k}^{g(\alpha)}
$$
As 
$$
\zeta_\nu < Z_{\nu+1}^{-\alpha},\,\,\,\,
\zeta_k < Z_{k+1}^{-\alpha}
$$
this will be enough to obtain Theorem 1.

Suppose that $g>\alpha -1$.
Either 
$$
Z_{k+1}\ge Z_k^g
$$
or
$$
Z_{k+1}< Z_k^g
,$$
The last inequality 
together with the inequality (\ref{ge1}) of  Lemma 5
gives
\begin{equation}\label{dop}
Z_k \ge Z_\nu^{\frac{\alpha}{g-\alpha +1}}.
\end{equation}
Now from  inequality (\ref{norma1}) of Corollary to Lemma 3
and (\ref{u2},\ref{dop}) we get
$$
\zeta_\nu\le
Z_\nu^{-\max (\alpha-1,1)-\frac{\alpha(\alpha-1)}{g-\alpha+1}}.
$$
So
$$
\beta (\Theta)\ge
\max_{ g>\alpha-1}
\,\,
\min \, \left(
\alpha g, \max (\alpha-1,1)+\frac{\alpha(\alpha-1)}{g-\alpha+1}
\right).
$$
But $g(\alpha )$ is the solution of (\ref{eq}).
So Theorem is proved.$\Box$

 \vskip+1.0cm

author:

\vskip+0.3cm

Nikolay  G. Moshchevitin

e-mail: moshchevitin@mech.math.msu.su, moshchevitin@gmail.com

\end{document}